\documentclass[final]{siamltex}

\usepackage[english]{babel}
\usepackage{amsmath,amssymb,amsfonts}
\usepackage[ruled]{algorithm2e} 
\usepackage{soul}
\usepackage{verbatim}
\usepackage{pgfplots}

\newtheorem{thm}{Theorem}

\newtheorem{remark}[thm]{\textit{Remark}}

\graphicspath{{figures/}}
\newcommand{\ee}{{\rm e}\hspace{1pt}}
\newcommand{\ii}{\text{i}\hspace{1pt}}
\newcommand{\dd}{\hspace{1pt}{\rm d}\hspace{0.5pt}}

\newcommand{\norm}[1]{\Vert #1 \Vert}
\newcommand{\normm}[1]{\left\lVert#1\right\rVert}

 \title{Analysis of Krylov subspace approximation to \\
Large Scale Differential Riccati Equations }

\author{Antti Koskela\thanks{Department of Mathematics and Statisctics, University of Helsinki, \texttt{antti.h.koskela@helsinki.fi}. }
\and Hermann Mena\thanks{Department of Mathematics, Yachay Tech, Urcuqu\'i, Ecuador, Department of Mathematics, University of Innsbruck, Innsbruck, Austria, \texttt{mena@yachaytech.edu.ec}.}
}

    \date{}

\begin{document}
   \maketitle

\begin{abstract}
We consider a Krylov subspace approximation method for the symmetric differential Riccati equation
$\dot{X} = AX + XA^T + Q - XSX$, $X(0)=X_0$. The method we consider is based on projecting
the large scale equation onto a Krylov subspace spanned by the matrix $A$ and the low rank factors of $X_0$ and $Q$. 
We prove that the method is structure preserving in the sense that it preserves 
two important properties of the exact flow, namely the positivity of the exact flow,
and also the property of monotonicity. 
We also provide a theoretical a priori error analysis which 
shows a superlinear convergence of the method. This behavior is illustrated in the numerical experiments. 
Moreover, we derive an efficient a posteriori error estimate
as well as discuss multiple time stepping combined with a cut of the rank of the numerical solution.
\end{abstract}



\section{Introduction}

 Large scale differential Riccati equations (DREs) arise in the numerical treatment of optimal 
 control problems governed by partial differential equations.
 This is the case in particular when solving a linear quadratic regulator problem (LQR), 
 a widely studied problem in control theory. 
 We shortly describe the finite dimensional LQR problem. For more details, we refer
 to~\cite{Freiling_book,Corless_book}. 
The differential Riccati equation arises in the finite horizon case, i.e., 
 when a finite time integral cost functional is considered. Denoting the time interval $[0,t_f]$,
 the functional has then the quadratic form
\begin{equation} \label{eq:quad_cost}
  J(x,u) = \int\limits_0^{t_f} \big(x(t)^T C^T C x(t) + u(t)^T u(t) \big) \, \dd t + x(t_f)^T G x(t_f), 
\end{equation}
where $x \in \mathbb{R}^n$, $C \in \mathbb{R}^{q \times n}$ ($q \ll n$) and $u \in \mathbb{R}^r$ ($r \ll n$).
The coefficient matrix $G$ of the penalizing term $x(t_f)^T G x(t_f)$ is symmetric, nonnegative and has a low rank.  
The functional \eqref{eq:quad_cost} is constrained by the system of differential equations
\begin{equation}  \label{eq:lin_op}
  \dot{x}(t) = A x(t) + B u(t), \qquad  x(0)= x_0, \quad  t\in[0,t_f]
\end{equation}
where the matrix $A \in \mathbb{R}^{n \times n}$ is sparse and $B\in \mathbb{R}^{n \times r}$. 
The number of columns of $B$ corresponds to the number of controls and the matrix $C$ represents an observation matrix. 
Under suitable conditions \cite{Freiling_book,Corless_book}, the control $\widetilde{u}$ 
minimizing the functional \eqref{eq:quad_cost} is given by
\begin{equation} \label{eq:control}
\widetilde{u}(t) = K(t) \widetilde{x}(t), \quad \textrm{where} \quad K(t) = -B^T X(t),
\end{equation}
$X(t)$ is the unique solution of 
\begin{equation} \label{eq:optimal_DRE}
\dot{X} + A^T X + X A - X B B^T X + C^T C = 0, \quad X(t_f) = G,
\end{equation}
and  $\widetilde{x}(t)$ satisfies
$$
\dot{\widetilde{x}}(t) = \big(A - B B^T X(t_f-t) \big) \widetilde{x}(t), \quad \widetilde{x}(0) = x_0.
$$
As a result, the central computational problem becomes that of solving the final value problem \eqref{eq:optimal_DRE} which, 
with a careful change of variables, can be written as a initial value problem.

We consider a Krylov subspace approximation method for large scale differential Riccati equations of the form \eqref{eq:optimal_DRE}.
A similar projection method for DREs has been recently proposed
in~\cite{Guldogan}. Our approach differs from that of~\cite{Guldogan} in the fact that the initial value matrix $G$
of \eqref{eq:optimal_DRE} is contained in the Krylov subspace. This allows multiple time stepping.
Our approach is also related to projection techniques considered for large scale algebraic Riccati  equations~\cite{Lin,Simoncini2016}.

Essentially, the method we consider is based on projecting the matrices $A,Q,S$ and $X_0$ on an appropriate Krylov subspace,
namely on the \it block Krylov subspace \rm spanned by $A$ and the low rank factors of $X_0$ and $Q$.
The projected small dimensional system is then solved using existing linearization techniques. 
We show that when using a Pad\'e approximant to solve the linearized
small dimensional system, the total approximation will be structure preserving in a sense 
that the property of the positivity is preserved. Also the property of monotonicity is preserved under certain conditions.
Our Krylov subspace approach is also strongly related to Krylov subspace techniques used for
approximation of the product of an matrix function and a vector, $f(A)b$,
and to exponential integrators~\cite{HochbruckOstermann}.
For an introduction to matrix functions we refer to the monograph~\cite{Higham}.
The effectiveness of these techniques comes from the fact that generating Krylov subspaces is essentially based 
on operations of the form $b \rightarrow Ab$, which are cheap for sparse $A$.

The linearization approach for DREs is a well-known method. 
This allows an efficient integration for dense problems, see e.g.~\cite{Laub}.
Another approach, the so called Davison--Maki method~\cite{DavisonMaki}, uses the fundamental solution of the linearized system. 
A modified variant, avoiding some numerical instabilities, is proposed in~\cite{KenneyLeipnik}. 
However, the application of these methods for large scale problems is impossible due to the
high dimensionality of the linearized differential equation.

The problem of solving large scale DREs has received recently considerable attention.
In~\cite{BennerMena1, BennerMena2} the authors proposed efficient BDF and Rosenbrock methods for solving DREs capable of exploiting
several of the above described properties: sparsity of $A$, low rank structure of $B$, $C$ and $G$, and the symmetry of the solution.
However, several difficulties arise when approximating the optimal control \eqref{eq:control}
in the large scale setting. One difficulty is to evaluate the state equation $x(t)$ and Riccati equation 
$X(t)$ in the same mesh. 
In~\cite{Saak2} a refined ADI integration method is proposed which addresses the high storage requirements
of large scale DRE integrators.
In recent studies an efficient splitting method~\cite{Stillfjord}  
and adaptive high-order splitting schemes~\cite{Stillfjord2} for large scale DREs have been proposed.

The paper is organized as follows. In Section 2 we describe some preliminaries. Then, in Section 3, the structure preserving method is proposed. 
In Section 4, the error analysis first for the differential Lyapunov equation (a simplified version of the DRE), and then for the DRE is presented. 
In Section 5 a posteriori error estimation is described. In Section 6 the rank cut and multiple time stepping are discussed. 
Numerical examples and conclusions of Sections 7 and 8 conclude the article.

\paragraph{Notation and definitions}
Throughout the article $\norm{ \cdot }$ will denote the Euclidean norm, or its induced matrix norm, i.e., the spectral norm.
By $\mathrm{R}(A)$ we denote the column space of a matrix $A$.
We say that a matrix $A$ is nonnegative if it is symmetric positive semidefinite, and write $A \geq 0$.
For symmetric matrices $A$ and $B$ we write $B \geq A$ if $B - A \geq 0$.

We will repeatedly use the notion of the \it logarithmic norm \rm of a matrix
$A \in \mathbb{C}^{n \times n}$. It can be defined via the \it field of values \rm $\mathcal{F}(A)$, which
is defined as
 $$
 \mathcal{F}(A) = \{ x^* A x \, : \, x \in \mathbb{C}^n, \, \norm{x} = 1 \}.
 $$
 Then, the logarithmic norm $\mu(A)$ of $A$ is defined by
 $$
 \mu(A) := \{ \max \, \textrm{Re} \, z \, : \, z \in \mathcal{F}(A) \}.
 $$
 We will also repeatedly use the exponential-like function $\varphi_1$ defined by
 $$
 \varphi_1(z) = \frac{\ee^z - 1}{z} = \sum\limits_{\ell=0}^\infty \frac{z^\ell}{(\ell + 1)!}.
 $$

\section{Preliminaries}
 
From now on we consider the time invariant symmetric differential Riccati equation (DRE) written in the form
\begin{equation} \label{eq:DRE}
\begin{aligned}
\dot{X}(t) &=  A X(t) + X(t) A^T + Q - X(t) S X(t),\\
X(0) &= X_0,
\end{aligned}
\end{equation}
where $t \geq 0$ and $A,Q,S,X_0 \in \mathbb{R}^{n \times n}$, $Q^T = Q$, $S^T = S$.
Specifically, we focus on the low rank positive semidefinite case, where 
\begin{equation} \label{eq:initial_data}
 X_0 = Z Z^T, \quad  Q = C C^T,
\end{equation}
for some $Z \in \mathbb{R}^{n \times p}$ and $C \in \mathbb{R}^{n \times q}$, $p,q \ll n$, and $S$ is positive semidefinite.
Notice that we changed here from to $A^T$ to $A$ (a common choice the numerical analysis literature~\cite{Dieci_Eirola,Dieci_Eirola2})
and from now on $C$ is tall and skinny instead of short and fat as in \eqref{eq:optimal_DRE}.
Although $S$ arises from the low rank matrix $B$ in \eqref{eq:optimal_DRE}, we do not place any restriction on the rank of $S$.

\subsection{Linearization}

We recall a fact that will be needed later on (see e.g.~\cite[Thm.\;3.1.1.]{Freiling_book}).
\begin{lemma}[Associated linear system] \label{lem:linearization}
The DRE \eqref{eq:DRE} is equivalent to solving the linear system of  differential equations
\begin{equation} \label{eq:hamdiff}
\frac{\dd}{\dd t}\begin{bmatrix}
U(t) \\ V(t)
\end{bmatrix} =
\begin{bmatrix}
-A & S \\ Q  & A^T
\end{bmatrix}
\begin{bmatrix}
U(t) \\ V(t)
\end{bmatrix}, \quad
\begin{bmatrix}
U(0) \\ V(0)
\end{bmatrix} =
\begin{bmatrix}
I \\ X_0
\end{bmatrix}
\end{equation}
 where $U(t),V(t) \in\mathbb{R}^{n\times n}$.
If the solution $X(t)$ of \eqref{eq:DRE} exists on the interval $[0, T]$, then the solution of \eqref{eq:hamdiff}
exists, $U(t)$ is invertible on $[0, T]$, and
\begin{equation*} 
X(t)=V(t)U(t)^{-1}.
\end{equation*}
\end{lemma}
Notice also that the matrix $\mathcal{H} = \begin{bmatrix}
-A & S \\ Q  & A^T
\end{bmatrix}$ is Hamiltonian, i.e., it holds that
\begin{equation} \label{eq:Ham_property}
(J \mathcal{H})^T = J \mathcal{H}, \quad \textrm{where} \quad J = \begin{bmatrix} 0  & I \\ -I & 0 \end{bmatrix}.
\end{equation}

This linearization approach 
is a standard method for 
solving finite dimensional DREs, and leads to efficient integration methods for dense problems, 
see e.g.~\cite{DavisonMaki}.

\subsection{Integral representation of the exact solution}

For the exact solution of \eqref{eq:DRE} we have the following integral representation (see also~\cite[Thm.\;8]{Kucera}).
\begin{theorem}[Exact solution of the DRE] \label{thm:exact_solution}
The exact solution of the DRE \eqref{eq:DRE} is given by
\begin{equation} \label{eq:exact_solution}
\begin{aligned}
X(t) = \ee^{t A} X_0 \ee^{t A^T} & + \int\limits_0^t \ee^{(t-s) A} Q \ee^{(t-s) A^T} \, \dd s\\ 
& - \int\limits_0^t \ee^{(t-s) A} X(s) S X(s) \ee^{(t-s) A^T} \, \dd s.
\end{aligned}
\end{equation}
\begin{proof} The proof can be carried out by elementary differentiation. \end{proof}
\end{theorem}

\subsection{Positivity and monotonicity of the exact flow}

We recall two important properties of the symmetric DRE, namely the positivity
 of the exact solution (see e.g.~\cite[Prop.\;1.1]{Dieci_Eirola})
and the monotonicity of the solution with relative to the initial data (see e.g.~\cite[Thm.\;2]{Dieci_Eirola2}). 
By these we mean the following.

\begin{theorem}[Positivity and monotonicity of the solution] \label{thm:positivity_and_monotonicity}
For the solution $X(t)$ of the symmetric DRE \eqref{eq:DRE} it holds:

\vspace{2mm}

\begin{enumerate}
 \item (Positivity) $X(t)$ is symmetric positive semidefinite and it exists for all $t>0$.

\item (Monotonicity) Consider two symmetric DREs of the \eqref{eq:DRE} corresponding 
to the linearized systems of the form \eqref{eq:hamdiff} with the coefficient matrices
$$
\mathcal{H}  = \begin{bmatrix} - A & S \\ Q & A^T \end{bmatrix} \quad \textrm{and} \quad \widetilde{\mathcal{H}}  
= \begin{bmatrix}  -\widetilde{A} & \widetilde{S} \\ \widetilde{Q} & \widetilde{A}^T \end{bmatrix}
$$
and let $J$ be the skew-symmetric matrix \eqref{eq:Ham_property}.
Then, if $\widetilde{\mathcal{H}} J \leq \mathcal{H} J$, and if $0 \leq X_0 \leq \widetilde{X}_0$,
then for every $t\geq 0$:  $X(t) \leq \widetilde{X}(t)$.
\end{enumerate}
\end{theorem}

We will later show that our proposed numerical method preserves the properties of Theorem~\ref{thm:positivity_and_monotonicity}.

\subsection{Bound for the solution}

Using the positivity property of $X(t)$ (Thm.~\ref{thm:positivity_and_monotonicity}) we obtain the following bound
for the norm of the solution. This will be repeatedly needed in the analysis of the proposed method.

\begin{lemma}[Bound for the exact solution] \label{lem:bound_exact}
For the solution $X(t)$ of \eqref{eq:DRE} it holds
\begin{equation} \label{eq:bound_exact}
\norm{X(t)} \leq \ee^{2 t \mu(A)} \norm{X_0} + t \varphi_1 \big(2 t \mu(A) ) \norm{Q}.
\end{equation}

\begin{proof} Since $X_0$, $Q$ and $X(t)$ are all symmetric positive semidefinite, we see that the first two terms on the right hand side
of \eqref{eq:exact_solution} are symmetric positive semidefinite and the third term is symmetric negative semidefinite.
Moreover, since $X(t)$ is symmetric positive semidefinite by Theorem~\ref{thm:positivity_and_monotonicity}, 
and since for every symmetric positive definite matrix $M$ it holds that $\norm{M} = \max\limits_{\norm{x}=1} x^* M x$,
we see that
\begin{equation*}
\begin{aligned}
\norm{X(t)} \leq \norm{\ee^{t A} X_0 \, \ee^{t A^T} + \int\limits_0^t \ee^{(t-s) A} Q \, \ee^{(t-s) A^T} \, \dd s}.
\end{aligned}
\end{equation*}
Using the well-known bound $\norm{\ee^{t A}} \leq \ee^{t \mu(A)}$ (see e.g.~\cite[p.\;138]{Embree_Trefethen}), 
the fact that $\mu(A^T) = \mu(A)$ and that
$t \varphi_1(t z) = \int_0^t \ee^{(t-s)z} \,\dd s$, the claim follows. \end{proof}
\end{lemma}

From Lemma~\ref{lem:bound_exact} we immediately get the following corollary.
\begin{corollary} \label{cor:max_X}
The solution $X(t)$ satisfies
$$
\max_{s \in [0,t]} \norm{X(s)} \leq \max\{1, \ee^{2 t \mu(A)} \} \norm{X_0} +t \max \{1, \varphi_1 \big(2 t \mu(A)) \} \norm{Q}.
$$
\end{corollary}

\section{A Krylov subspace approximation and its structure preserving properties}

In this section we propose our projection method. The original problem \eqref{eq:DRE} is projected to
small dimensional space using a matrix $V_k$ with orthonormal columns which contains certain Krylov subspaces.
The fact that $V_k$ needs to contain these subspaces
can be seen from the point of view of Krylov subspace approximation of the matrix exponential
(see also the solution formula \eqref{eq:exact_solution}).
This is strongly related to the approach taken by Saad already in~\cite{Saad_Lyapunov}
for the algebraic Lyapunov equation.
Before introducing our projection method, we recall some basic facts about the Krylov subspace approximation of the matrix exponential.
This will also give some auxiliary results that are needed later in the convergence analysis.

\subsection{Block Krylov subspace approximation of the matrix exponential} \label{subsec:Krylov_approx}

The Krylov subspace approximation of products of the form $f(A)b$ has recently been 
an active topic of research, and we mention the work on classical Krylov subspaces 
\cite{DruskinKnizhnerman,GallopoulosSaad,Knizhnerman,Saad}, 
extended Krylov subspaces \cite{Knizhnerman}, and rational Krylov subspaces \cite{EshofHochbruck,Beckermann_Reichel}.

Block Krylov subspace methods are based on the idea of projecting a high dimensional
problem involving a matrix $A \in \mathbb{R}^{n \times n}$ and a block matrix
$B \in \mathbb{R}^{n \times \ell}$ onto a lower dimensional subspace, a block Krylov subspace $\mathcal{K}_k (A,B)$,
which is defined by
\begin{equation} \label{def:krylov}
\mathcal{K}_k(A,B) = \textrm{span} \{B, AB, A^2 B, \ldots, A^{k-1} B \}.
\end{equation}
Usually, an orthogonal basis matrix $V_k$ for $\mathcal{K}_k(A,B)$ is
generated using an Arnoldi type iteration, and this matrix is then used for the projections.
There exist several Arnoldi type methods to produce an orthogonal basis matrix for $\mathcal{K}_k(A,B)$,
and in numerical experiments we use the \it block Arnoldi iteration \rm given in~\cite{Saad_book} which is listed algorithmically as follows.

\bigskip

\begin{enumerate}

 \item \textbf{Input:} $A \in \mathbb{R}^{n \times n}$, $B \in \mathbb{R}^{n \times \ell}$ and number of iterations $k$.

 \item \textbf{Start.} Compute QR decomposition: $B = U_1 R_1$.
  
 \item \textbf{Iterate.} \it for $j=1,...,k$ compute: \\
 \begin{equation*}
 \begin{aligned}
 H_{ij} &= U_i^T A U_j, \quad i = 1, \ldots, j, \\
W_j &= A U_j - \sum\limits_{i=1}^j U_i H_{ij},  \\
W_j &= U_{j+1} H_{j+1,j}\quad (\textrm{QR decomposition of } W_j ). 
 \end{aligned}
 \end{equation*}
 \end{enumerate} 

\medskip

As usual, the orthogonalisation can be carried out at step 3 in a modified Gram--Schmidt manner and reorthogonalisation
can be performed if needed. 

This algorithm gives a basis matrix with orthogonal columns, $V_k = \begin{bmatrix}
U_1 & \ldots & U_k \end{bmatrix} \in \mathbb{R}^{n \times k \ell}$, for the block Krylov subspace
$\mathcal{K}_k(A,B)$ and the projected block Hessenberg matrix 
\begin{equation} \label{eq:block_Hessenberg}
H_k = V_k^T A V_k.
\end{equation}
This means that the $\ell \times \ell$ $(i,j)$-block of $H_k$ is given by $H_{ij}$ in the above algorithm.
Moreover, the following Arnoldi relation holds:
\begin{equation}\label{eq:lanczos-it}
A V_k = V_k H_k + U_{k+1} H_{k+1,k}  E_k^T,
\end{equation}
where $E_k = \begin{bmatrix} 0 & \ldots & 0 & I_\ell \end{bmatrix}^T \in \mathbb{R}^{k \ell \times \ell}$.

If $A$ has its field of values on a line, e.g., is Hermitian or skew-Hermitian, 
then there exists $\theta \in \mathbb{R}$ such that $\ee^{\ii \theta} A$ is Hermitian.
By \eqref{eq:block_Hessenberg} this implies that $H_k$ is block tridiagonal, the orthogonalisation recursions become three-term recursions,
and we get the \it block Lanczos iteration. \rm

%
%


The polynomial approximation property of Krylov subspaces motivates to approximate 
the product of the matrix exponential and a block matrix as
\begin{equation} \label{eq:exp_approx}
\ee^A B \approx V_k \ee^{H_k} V_k^T B =  V_k \ee^{H_k} E_1 R_1,
\end{equation}
where $E_1 = \begin{bmatrix} I_\ell & 0 & \ldots & 0 \end{bmatrix}^T \in \mathbb{R}^{ k \ell \times \ell  }$. 
For a vector $B$, the approximation \eqref{eq:exp_approx} was
considered already in~\cite{DruskinKnizhnerman,GallopoulosSaad},
and for the case of a block matrix $B$ it has been considered also in~\cite{Lopez_Simoncini}.

Since the columns of $V_k$ are orthonormal, we have
$
\mathcal{F}(H_k) = \mathcal{F}(V_k^T A V_k)\subset \mathcal{F}(A).
$
and from this it follows that $\mu(H_k) \leq \mu(A)$. 
Clearly, it also holds that $\norm{H_k} \leq \norm{A}$. 
Moreover, we have the following bound. 
\begin{lemma} \label{lem:exp_error}
For the approximation \eqref{eq:exp_approx} holds
\begin{equation} \label{eq:exp_error}
\norm{ \ee^{t A } B - V_k \ee^{t H_k } V_k^T B } \leq 2 \max\{1,\ee^{t \mu(A)}\} \frac{\norm{t A}^k}{k ! }\norm{B}.
\end{equation}
\begin{proof}  The proof goes analogously to the proof of~\cite[Thm\;2.1]{GallopoulosSaad}, where $B$ is a vector. \end{proof}
\end{lemma}


\subsection{Rational Krylov subspaces}

We also mention the possibility of approximating matrix functions in \it rational Krylov subspaces \rm 
(see e.g.~\cite{DruskinLiebermann}, ~\cite{Guettel}, ~\cite{EshofHochbruck} and ~\cite{Simoncini2016}).
For poles $\bar{s} = \{s_1, s_2, \ldots \}$, $s_i \in \mathbb{C}$, the rational Krylov subspace can be defined
as (see also)
\begin{equation} \label{def:rational}
\mathcal{K}_k(A,B,\bar{s}) = \textrm{span} \{B, (s_1 I- A)^{-1} B,  \ldots,  \prod\limits_{\ell=1}^{k-1} (s_\ell I- A)^{-1} B \}.
\end{equation}
Then, if a matrix $V_k$ with orthogonal columns gives a basis for the subspace $\mathcal{K}_k(A,B,\bar{s})$,
the matrix exponential can be approximated as \eqref{eq:exp_approx}, where $H_k = V_k^T A V_k$.
 Especially for sparse matrices,
the rational Krylov methods are often more efficient, and as the solution usually converges faster
with respect to subspace dimension, the rational alternative is usually more memory efficient.
These differences will be illustrated in numerical experiments.
However, for simplicity, in our analysis and numerical experiments  we will use  the block Arnoldi iteration.

\subsection{The method}

We approximate $X(t)$ in 
the block Krylov subspace $\mathcal{K}_k \big(A, \begin{bmatrix} Z & C \end{bmatrix} \big)$.
The fact that the projection onto this subspace results as an accurate approximation can be seen from the 
form of the exact solution \eqref{eq:exact_solution} and from the Krylov approximation properties shown in the
last subsection.
To this end, an orthogonal basis matrix $V_k \in \mathbb{R}^{n \times k (p+q)}$ for $\mathcal{K}_k \big(A, \begin{bmatrix} Z & C \end{bmatrix} \big)$ 
is first generated using the block Arnoldi iteration.
Then, we carry out the approximation as listed in Algorithm~\ref{Alg:main}.
Notice that the method works independently of the rank of $S$.

\begin{algorithm} \label{Alg:main} 
\caption{Krylov subspace approximation of the DRE \eqref{eq:DRE}}
\SetKwInOut{Input}{Input}\SetKwInOut{Output}{Output}

\vspace{2mm}

\Input{Time step size $h$, Krylov subspace size $k$, matrices $A,S \in \mathbb{R}^{n \times n}$, 
$Z \in \mathbb{R}^{n \times p}$ and $C \in \mathbb{R}^{n \times q}$, $p,q \ll n$.}

\nl Carry out $k$ steps of the block Arnoldi iteration to obtain

  \begin{itemize}
   \item the orthogonal basis matrix $V_k$ of $\mathcal{K}_k \big(A, \begin{bmatrix} Z & C \end{bmatrix} \big)$
   \item the block Hessenberg matrix  $H_k = V_k^T A V_k$
   \item the matrices $C_k = V_k^T C$ and $Z_k = V_k^T Z$
  \end{itemize}
    
\nl Compute $S_k = V_k^T S V_k$.

\nl Compute the solution $Y_k(t)$ of the small dimensional system
\begin{equation} \label{eq:riccati_small_system}
\begin{array}{rcl}
\dot{Y}_k(t) &=&  H_k Y_k(t) + Y_k(t) H_k^T + C_k C_k^T - Y_k(t) S_k Y_k(t),\\
Y_k(0) &=& Z_k Z_k^T.
\end{array}
\end{equation}

\BlankLine

\nl Approximate 
$
X(t) \approx X_k(t) = V_k Y_k(t) V_k^T.
$

\end{algorithm}

\subsection{Solving the small dimensional system} \label{subsec:scaling_and_squaring}

To solve the small dimensional system \eqref{eq:riccati_small_system} 
we use the \it modified Davison--Maki method\rm~\cite{KenneyLeipnik}. This method is chosen because of its structure preservation properties 
 which are shown in Subsection~\ref{subsec:structure}. The method can be described as follows. 

As shown in Lemma~\ref{lem:linearization}, the solution of the projected system \eqref{eq:riccati_small_system}
is given by
\begin{equation} \label{eq:exp_H_k}
Y_k(t) = W_k(t) U_k(t)^{-1}, \textrm{ where} \quad \begin{bmatrix}
 U_k(t) \\ W_k(t)
\end{bmatrix} = \exp\left( t \begin{bmatrix}
              - H_k & S_k  \\  C_k C_k^T &  H_k^T
             \end{bmatrix} \right)
\begin{bmatrix}
I_k \\ Z_k Z_k^T 
\end{bmatrix}.
\end{equation}
Instead of directly evaluating 
$Y_k(t)$ by \eqref{eq:exp_H_k},
which is the idea of the original Davison--Maki method~\cite{DavisonMaki}, we perform substepping in order to avoid
numerical instabilities arising from the inversion of the matrix $U_k(t)$ in \eqref{eq:exp_H_k}.
This is exactly the modified Davison--Maki method, and it is presented in the following pseudocode.
We denote  $Y_k^j \approx Y_k \big(\tfrac{j \cdot t}{m} \big)$.

\bigskip

\begin{enumerate}

 \item \textbf{Input:} Hamiltonian matrix $\left[ \begin{smallmatrix}
               - H_k & S_k  \\  C_k C_k^T &  H_k^T
             \end{smallmatrix} \right]$, $Y_k(0) = Z_k Z_k^T$, \\ 
             time $t>0$, substep size $\Delta t = t/m$,  $m \in \mathbb{Z}_+$. \\
  
 \item \textbf{Set:} $Y_k^0 = Y_k(0)$.  \\
  
 \item \textbf{Iterate.} \it for $j=0,...,m-1$: \\
$$
\begin{bmatrix}
 U_{j+1}  \\ W_{j+1} 
\end{bmatrix}
=
\exp\left( \Delta t \begin{bmatrix}
              - H_k & S_k  \\  C_k C_k^T &  H_k^T
             \end{bmatrix} \right)
\begin{bmatrix}
  I_k   \\   Y_k^j 
\end{bmatrix}, \quad \quad Y_k^{j+1} = W_{j+1}  U_{j+1}^{-1}.
$$

 \end{enumerate} 

\medskip

For computing the matrix exponential in Step 3, we use the 13$th$ order diagonal Pad\'e aproximant
which is implemented in Matlab as 'expm' command~\cite{Higham2}. \\

\subsection{Structure preserving properties of the approximation} \label{subsec:structure}

We next inspect the two properties stated in Theorem~\ref{thm:positivity_and_monotonicity}.
We show that the proposed projection method preserves the property of the positivity of the exact flow, 
and it also preserves the property of monotonicity under the condition that the matrix $V_k$ used
for the projection stays constant when the initial data for the DRE is changed.
Notice that these results are not restricted to polynomial Krylov subspace methods.

\begin{theorem} \label{thm:numerical_positivity}
The numerical approximation given by Algorithm~\ref{Alg:main} preserves the property of positivity
stated in Theorem~\ref{thm:positivity_and_monotonicity}.
\begin{proof} The projected coefficient matrices $S_k$, $C_kC_k^T$ and the initial value $Z_k Z_k^T$ of the small system
\eqref{eq:riccati_small_system} are clearly all symmetric nonnegative.
Thus the small system \eqref{eq:riccati_small_system} is a symmetric DRE.
By Theorem 3.1 of~\cite{Dieci_Eirola}, an application of a symplectic Runge--Kutta scheme with positive weights
$b_i$ (see~\cite{Dieci_Eirola2} for details) gives as a result a symmetric nonnegative solution. As the $s$th order diagonal Pad\'e approximant
equals the stability function of the $s$-stage Gauss--Legendre method (see e.g.~\cite[p.\;46]{Iserles}),
the Pad\'e approximation in the third substep of the modified Davison--Maki method (Subsection~\ref{subsec:scaling_and_squaring})
corresponds to a symplectic Runge--Kutta method.
Thus each substep of the modified Davison--Maki method
outputs a symmetric nonnegative matrix and as a result $Y_k(t)$ is symmetric nonnegative. 
Therefore also $X_k(t) = V_k Y_k(t) V_k^T$ is symmetric nonnegative. 
\end{proof}
\end{theorem}

\begin{theorem}
The numerical approximation given by Algorithm~\ref{Alg:main} preserves the property of monotonicity
in the following sense. 
Consider two DREs corresponding to linearizations with the coefficient matrices
\begin{equation} \label{eq:linearizations}
\mathcal{H}  = \begin{bmatrix} - A & S \\ Q & A^T \end{bmatrix} \quad \textrm{and} \quad \widetilde{\mathcal{H}}  
= \begin{bmatrix} - \widetilde{A} & \widetilde{S} \\ \widetilde{Q} & \widetilde{A}^T \end{bmatrix}
\end{equation}
such that 
\begin{equation} \label{eq:Ham_ineq}
\widetilde{\mathcal{H}} J \leq \mathcal{H} J, \quad 0 \leq X_0 \leq \widetilde{X}_0.
\end{equation}
Suppose both systems are projected using the same orthogonal matrix $V_k \in \mathbb{R}^{n \times k}$,
giving as a result small $k$-dimensional systems of the form \eqref{eq:riccati_small_system}
for the matrices $Y_k(t)$ and $\widetilde{Y}_k(t)$. Then, for the matrices $X_k(t) = V_k Y_k(t) V_k^T$
and $\widetilde{X}_k(t) = V_k \widetilde{Y}_k(t) V_k^T$ we have
$$
X_k(t) \leq \widetilde{X}_k(t).
$$
\begin{proof} Consider the projected systems of the form \eqref{eq:riccati_small_system} corresponding to $Y_k(t)$ and $\widetilde{Y}_k(t)$
with the projected coefficient matrices $H_k$, $Q_k$ and $S_k$, and  $\widetilde{H}_k$, $\widetilde{Q}_k$ and $\widetilde{S}_k$,
respectively. Consider also the corresponding linearizations of the form \eqref{eq:linearizations} 
with the Hamiltonian matrices 
$$
\mathcal{H}_k :=  \begin{bmatrix} V_k & 0 \\ 0& V_k \end{bmatrix}^T \mathcal{H} \begin{bmatrix} V_k & 0 \\ 0& V_k \end{bmatrix}
\quad \textrm{and} \quad  \widetilde{\mathcal{H}}_k := \begin{bmatrix} V_k & 0 \\ 0& V_k \end{bmatrix}^T \widetilde{\mathcal{H}} \begin{bmatrix} V_k & 0 \\ 0& V_k \end{bmatrix}.
$$
By the reasoning of the proof of Theorem~\ref{thm:numerical_positivity}, the projected systems corresponding
to $Y_k(t)$ and $\widetilde{Y}_k(t)$ are symmetric DREs.
We see that
$$
\widetilde{\mathcal{H}}_k J_k - \mathcal{H}_k J_k =  \begin{bmatrix} V_k & 0 \\ 0& V_k \end{bmatrix}^T
( \widetilde{\mathcal{H}} J - \mathcal{H} J )  \begin{bmatrix} V_k & 0 \\ 0& V_k \end{bmatrix},
$$
where $J_k = \left[ \begin{smallmatrix} 0 & I \\ -I & 0 \end{smallmatrix} \right] \in \mathbb{R}^{2k \times 2k}$.
Thus, from \eqref{eq:Ham_ineq} it follows that $\widetilde{\mathcal{H}}_k J_k \leq \widetilde{\mathcal{H}}_k J_k$.
Clearly, also $0 \leq Y_k(0) \leq \widetilde{Y}_k(0)$.
By Theorem 6 of \cite{Dieci_Eirola2}, an application of a symplectic Runge--Kutta scheme with positive weights
$b_i$ (see~\cite{Dieci_Eirola2} for details) preserves the monotonicity. 
Thus the Pad\'e approximants of
the substeps of the modified Davison--Maki method (Subsection~\ref{subsec:scaling_and_squaring})
preserve the monotonicity. Therefore, $Y_k(t) \leq \widetilde{Y}_k(t)$ and as a consequence $X_k(t) \leq \widetilde{X}_k(t)$. \end{proof}
\end{theorem}

\begin{remark}
As the basis matrix $V_k$ given by Algorithm~\ref{Alg:main} is independent of the
matrix $S=BB^T$ in the DRE \eqref{eq:DRE}, where $B$ is the control matrix in the original
linear system \eqref{eq:lin_op}, we see that Algorithm~\ref{Alg:main} preserves monotonicity
under modifications of $B$.  
However, if we change the initial value $X_0$ or the matrix $Q$, then forming a new basis $V_k$ is generally needed.
The fact that $V_k$ is independent of $B$ can also be seen by considering similar projection methods for the algebraic Riccati equation,
see e.g.~\cite{Simoncini_SIAM} and the references therein.
\end{remark}

\section{A priori error analysis}

We first consider the approximation of the DRE without the quadratic term $-X S X$, i.e.,
we consider the differential Lyapunov equation. This clarifies the presentation as the derived bounds will
be needed when we consider the approximation of the differential Riccati equation.
We note, however, that the bounds for the Lyapunov equation are applicable outside of the scope
of the optimal control problems, e.g., when considering time integration of an inhomogeneous matrix differential equation.

\subsection{Error analysis for the Lyapunov equation}

Consider the symmetric Lyapunov differential equation with
low rank initial data and constant low rank inhomogeneity,
\begin{equation} \label{eq:Lyapunov}
\begin{aligned}
\dot{X}(t) &=  A X(t) + X(t) A^T + C C^T,\\
X(0) &= Z Z^T,
\end{aligned}
\end{equation}
where $Z \in \mathbb{R}^{n \times p}$ and $C \in \mathbb{R}^{n \times q}$, $p,q \ll n$.
Then, the approximation is given by $X_k(t) = V_k Y_k(t) V_k^T$, where $Y_k(t)$ is a solution of the projected
system \eqref{eq:riccati_small_system} with $S=0$. For the error of this approximation we obtain the following bound.
\begin{theorem} \label{thm:Lyapunov}
Let $A \in \mathbb{R}^{n \times n}$, $Z \in \mathbb{R}^{n \times p}$, $C \in \mathbb{R}^{n \times q}$,
and let $X(t)$ be the solution of \eqref{eq:Lyapunov}.
Let $V_k \in \mathbb{R}^{n \times m(q+p)}$ be an orthogonal basis of the block Krylov subspace 
$\mathcal{K}_k \big(A, \begin{bmatrix} Z & C \end{bmatrix} \big)$.
Let $Y_k(t)$ be the solution of the projected system
\eqref{eq:riccati_small_system} with $S=0$, and let $X_k(t) = V_k Y_k(t) V_k^T$. Then,
$$
\norm{X(t) - X_k(t)}  \leq  4  \max \{1, \ee^{2 t \mu(A)} \} \norm{A}^k \left( \frac{t^k}{k!}\norm{X_0} + \frac{t^{k+1}}{(k+1)!} \norm{Q} \right).
$$
\begin{proof} Using the integral representation of Theorem~\ref{thm:exact_solution} for both $X(t)$ and $Y_k(t)$, we see that
$$
X(t) - X_k(t) = Err_{1,k}(t) + Err_{2,k}(t),
$$
where
\begin{equation} \label{eq:err1}
Err_{1,k}(t) = \ee^{t A} Z Z^T \ee^{t A^T} - V_k \ee^{ t H_k} V_k^T Z Z ^T V_k \ee^{ t H_k^T} V_k^T
\end{equation}
and
\begin{equation} \label{eq:err2}
\begin{aligned}
Err_{2,k}(t) = & \int\limits_0^t \ee^{(t-s) A} C C^T \ee^{(t-s) A^T} \, \dd s \\
 & - \int\limits_0^t V_k \ee^{(t-s) H_k} V_k^T C C^T  V_k \ee^{(t-s) H_k^T} V_k^T \, \dd s.
\end{aligned}
\end{equation}
Adding and substracting $ \ee^{t A}  Z Z^T V_k \ee^{ t H_k^T} V_k^T $ to the right hand side of \eqref{eq:err1} gives
\begin{equation*}
\begin{aligned}
Err_1(t) = & \,\, \ee^{t A} Z \left(\ee^{t A}  Z -  V_k \ee^{ t H_k} V_k^T Z   \right)^T \\
& -  \left( V_k \ee^{ t H_k} V_k^T Z  - \ee^{t A}  Z  \right) Z ^T V_k \ee^{ t H_k^T} V_k^T.
\end{aligned}
\end{equation*}
Using Lemma~\ref{lem:exp_error} to bound the norm of $\ee^{ t A } Z - V_k \ee^{ t H_k } V_k^T Z$,
and using the fact that $\mu(H_k) \leq \mu(A)$ and that $\norm{X_0} = \norm{Z Z^T} = \norm{Z}^2$, gives
$$
\norm{Err_1(t)} \leq 
  4 \, \max(1,\ee^{2 t \mu(A)}) \frac{\norm{t A}^k}{k ! } \norm{X_0}.
$$
Then, similarly, adding and substracting the term $\int_0^t \ee^{(t-s) A}  C C^T V_k \ee^{ (t-s) H_k^T} V_k^T \, \dd s$
to \eqref{eq:err2} and applying Lemma~\ref{lem:exp_error} shows that
\begin{equation*}
 \begin{aligned}
\norm{Err_2(t)} & \leq 
  4 \norm{Q}  \int\limits_0^t \max\{1,\ee^{2 (t-s) \mu(A)}\} \frac{\norm{(t-s) A}^k}{k ! } \, \dd s \\
  & \leq 4 \norm{Q} \max \{ 1,\ee^{2 t \mu(A)} \}  \norm{A}^k \frac{t^{k+1}}{(k+1)!}
 \end{aligned}
\end{equation*}
which shows the claim. \end{proof}
\end{theorem}

We note that the error bound of Theorem~\ref{thm:Lyapunov}, similarly to the bounds given in~\cite{GallopoulosSaad}, exhibits a hump before it starts to decrease
in case $\norm{tA} > 1$.
Improved bounds for special cases of the matrix $A$ are possible by using, e.g., results of~\cite{HochbruckLubich}.

\subsection{Refined error bounds for the Lyapunov equation} \label{subsec:refined}

Although Theorem~\ref{thm:Lyapunov} shows the superlinear convergence speed for the approximation of the Lyapunov equation \eqref{eq:Lyapunov}, 
sharper bounds can be obtained, e.g., by using the bounds given in~\cite{HochbruckLubich}.
As an example we consider the following. If $A$ is symmetric negative semi-definite with its spectrum inside the interval
$[-4\rho, 0]$, and $V_k$ is an orthonormal basis matrix for the block Krylov subspace
$\mathcal{K}_k(A,B)$, we have (see \cite[Thm.\;2]{HochbruckLubich}) for the error 
$\varepsilon_k := \norm{\ee^{tA} B - V_k \ee^{t H_k} V_k^T B}$ the bound
\begin{equation} \label{eq:HL_bounds}
\begin{aligned}
\varepsilon_k & \leq 10 \, \ee^{-k^2/(5\rho \, t)} \norm{B}, \quad \quad \sqrt{4 \rho \, t} \leq k \leq 2 \rho \, t, \\
\varepsilon_k & \leq 10 \, (\rho \, t)^{-1} \ee^{- \rho \, t} \left( \frac{\ee \rho \, t}{k} \right)^k \quad \quad \, k \geq 2 \rho \, t.
\end{aligned}
\end{equation}

Using \eqref{eq:HL_bounds} and following the proof of Theorem~\ref{thm:Lyapunov},
we get the following bound for the case of a symmetric negative semidefinite $A$.

\begin{theorem}
Let $A \in \mathbb{R}^{n \times n}$, $Z \in \mathbb{R}^{n \times p}$ and $C \in \mathbb{R}^{n \times q}$ define the Lyapunov equation \ref{eq:Lyapunov}.
Let $V_k \in \mathbb{R}^{n \times m(q+p)}$ be an orthogonal basis matrix of the subspace $\mathcal{K}_k(A, \begin{bmatrix} Z & C \end{bmatrix}$.
Let $Y_k(t)$ be the solution of the projected (using $V_k$) system
\eqref{eq:riccati_small_system} with $S=0$, and let $X_k(t) = V_k Y_k(t) V_k^T$. Then,
for the  error $\varepsilon_k := \norm{X(t) - X_k(t)}$ it holds that
\begin{equation} \label{eq:symm_refined_bound}
\begin{aligned}
\varepsilon_k  & \leq 20 \, \ee^{-k^2/(5\rho \, t)} \big( \norm{X_0} + t \norm{Q} \big), \quad \quad \quad \sqrt{4 \rho \, t} \leq k \leq 2 \rho \, t, \\
\varepsilon_k  &\leq 20 \, (\rho \, t)^{-1} \ee^{- \rho \, t} \left( \frac{\ee \rho \, t}{k} \right)^k \big( \norm{X_0} + t \norm{Q} \big),
\quad \quad \, k \geq 2 \rho \, t,
\end{aligned}
\end{equation}
\end{theorem}

\begin{figure}[h!]
 \begin{center}
 \includegraphics[scale=0.61]{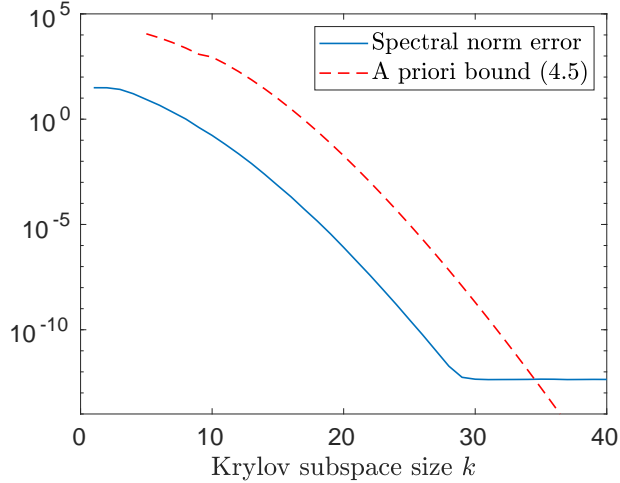}
 \end{center}
 \caption{Convergence of the approximation vs. the a priori bound given in the equation \eqref{eq:symm_refined_bound}. }
 \label{fig:bound1}
 \end{figure}

The bound \eqref{eq:symm_refined_bound} can be illustrated with the following simple numerical example.
Let $A \in \mathbb{R}^{400 \times 400}$ be the tridiagonal matrix $10^2 \cdot \mathrm{diag}(1,-2,1)$, $t=0.05$, and let
$Z\in \mathbb{R}^{400}$ and $C \in \mathbb{R}^{400}$ be random vectors. Figure \ref{fig:bound1} shows the convergence
of the algorithm vs. the a priori bound \eqref{eq:symm_refined_bound}.

\subsection{Error for the approximation of the Riccati equation} 
 
Here, we state our main theorem which shows the superlinear convergence property of Algorithm~\ref{Alg:main}
when applied to the DRE \eqref{eq:DRE}. Its proof, which is essentially based on Lemma~\ref{lem:exp_error} and Gr\"onwall's lemma,
is lengthy and is left to the appendix.

First, however, we state a bound for the norm of the numerical solution
$X_k(t)$ which will be needed in the proof of the main theorem.

\begin{lemma} \label{lem:bound_X_k}
Suppose, $X_0 = Z Z^T$, $Q = C C^T$ and that $S$ is symmetric nonnegative.
Then, $X_k(t)$ is symmetric nonnegative, and satisfies the bound
$$
\norm{X_k(t)} \leq \ee^{2 t \mu(A)} \norm{X_0} + t \varphi_1 \big(2 t \mu(A) ) \norm{Q}.
$$
\begin{proof} As $Z Z^T$, $C C^T$ and $S$ are symmetric and nonnegative, we see from \eqref{eq:riccati_small_system} that so are the
orthogonally projected matrices $Z_k Z_k^T$, $C_k C_k^T$ and $S_k$. Thus, the projected system is a symmetric DRE.
Applying Lemma~\ref{lem:bound_exact} to the projected system, and using the bounds $\mu(H_k) \leq \mu(A)$, $\norm{Q_k} \leq \norm{Q}$
and $\norm{V_k V_k^T X_0 V_k V_k^T} \leq \norm{X_0}$ shows the claim. \end{proof}
\end{lemma}

From Lemma~\ref{lem:bound_X_k} we immediately get the following bound.
\begin{corollary} \label{cor:max_X_k}
The numerical solution $X_k(t)$ satisfies
$$
\max_{s \in [0,t]} \norm{X_k(s)} \leq \max\{1, \ee^{2 t \mu(A)} \} \norm{X_0} +t \max \{1, \varphi_1 \big(2 t \mu(A)) \} \norm{Q}.
$$
\end{corollary}

We are now ready to state an error bound for the DRE. The proof is left to the appendix.

\begin{theorem} \label{thm:main_conv}
Let $A \in \mathbb{R}^{n \times n}$, $Z \in \mathbb{R}^{n \times p}$, $C \in \mathbb{R}^{n \times q}$ and $S \in \mathbb{R}^{n \times n}$
defined the DRE \eqref{eq:DRE}. Let $X_k(t)$ be the numerical solution given by Algorithm~\ref{Alg:main}. Then, the following bound holds:
\begin{equation} \label{eq:main_thm}
\norm{X(t) - X_k(t)}  \leq c(t) \norm{A}^k \left( \frac{t^k}{k!} \norm{X_0} + \frac{t^{k+1}}{(k+1)!}\norm{Q} \right),
\end{equation}
where
$$
c(t) = 4(1+2\norm{S} \alpha(t) \max \{ 1,\ee^{t \mu(A)}\} c_2(t)) \ee^{t \norm{S} \alpha(t)},
$$
$$
c_2(t) = 1 + t\norm{S} \alpha(t) \varphi_1\big( t \norm{S} \alpha(t) \max \{ 1,\ee^{t \mu(A)}\} \big)
$$
and
$$
\alpha(t) = \max\{1, \ee^{2 t \mu(A)} \} \norm{X_0} +t \max \{1, \varphi_1 \big(2 t \mu(A)) \} \norm{Q}.
$$
\end{theorem}

\section{A posteriori error estimation}

We consider next an a posteriori error estimation for the method by using ideas presented in~\cite{Botchev}. 

Denote the original DRE \eqref{eq:DRE} as
$$
\dot{X}(t) = F(X(t)), \quad X(0) = X_0.
$$
Using the residual matrix $R_k(t) = F(X_k(t)) - \dot{X}_k(t)$ we derive computable error estimates.
These derivations are based on the following lemma.	
\begin{lemma} \label{lem:apost_representation}
The error $\mathcal{E}_k(t) := X(t) - X_k(t)$ satisfies the equation
\begin{equation} \label{eq:apost_representation}
\begin{aligned}
\mathcal{E}_k(t) &= \int\limits_0^t \ee^{(t-s)A} R_k(s) \ee^{(t-s)A^T} \, \dd s  \\
& \quad \quad \quad - \int\limits_0^t \ee^{(t-s)A} \Big( \mathcal{E}_k(s) S X(s) + X_k(s) S \mathcal{E}_k(s) \Big) \ee^{(t-s)A^T} \, \dd s,
\end{aligned}
\end{equation}
where
\begin{equation} \label{eq:R_k_representation}
\begin{aligned}
R_k(t) = U_{k+1} H_{k+1,k} E_k^T Y_k(t) V_k^T  + V_k Y_k(t) E_k H_{k+1,k}^T U_{k+1}^T.
\end{aligned}
\end{equation}
\begin{proof} We see that the error $\mathcal{E}_k(t)$ satisfies the ODE
\begin{equation} \label{eq:E_k_ode}
\begin{aligned}
\dot{\mathcal{E}}_k(t) &= \dot{X}(t) - \dot{X}_k(t) = F(X(t)) - F(X_k(t)) + R_k(t) \\
&= A \big(X(t)-X_k(t)\big) + \big(X(t)-X_k(t)\big) A^T \\
& \quad \quad - X(t)SX(t) + X_k(t) S X_k(t) + R_k(t) \\
&= A \big(X(t)-X_k(t)\big) + \big(X(t)-X_k(t)\big) A^T \\
& \quad \quad - \big(X(t)-X_k(t)\big) S X(t) - X_k(t) S \big(X(t)-X_k(t)\big) + R_k(t) \\
&= A \mathcal{E}_k(t) + \mathcal{E}_k(t) A^T - \mathcal{E}_k(t) S X(t) - X_k(t) S \mathcal{E}_k(t) + R_k(t)
\end{aligned}
\end{equation}
with the initial value $\mathcal{E}_k(0) = 0$. Applying the variation-of-constants formula to \eqref{eq:E_k_ode} gives
\eqref{eq:apost_representation}.

Next we show the representation \eqref{eq:R_k_representation}. Since
$$
F(X_k(t)) = A V_k Y_k(t) V_k^T + V_k Y_k(t) V_k^T A^T + Q - V_k Y_k(t) V_k^T S V_k Y_k(t) V_k^T
$$
and 
$$
\dot{X}_k(t) = V_k H_k Y_k(t) V_k^T + V_k Y_k(t) H_k^T V_k^T +  V_k Q_k V_k^T - V_k Y_k(t) V_k^T S V_k Y_k(t) V_k^T
$$
we see that
\begin{equation} \label{eq:apost_med}
\begin{aligned}
R_k(t) & = F(X_k(t)) - \dot{X}_k(t) \\
& = (A V_k - V_k H_k) Y_k(t) V_k^T  + V_k Y_k(t) (A V_k - V_k H_k)^T + Q - V_k Q_k V_k^T \\
& = (A V_k - V_k H_k) Y_k(t) V_k^T  + V_k Y_k(t) (A V_k - V_k H_k)^T,
\end{aligned}
\end{equation}
since $V_k Q_k V_k^T = V_k V_k^T C C^T V_k V_k^T = C C^T = Q$ as $C \in \mathrm{R}(V_k)$.
Substituting the Arnoldi relation $A V_k - V_k H_k = U_{k+1} H_{k+1,k} E_k^T$ 
into \eqref{eq:apost_med} gives the representation \eqref{eq:R_k_representation}. \end{proof}
\end{lemma}

\begin{figure}[h!]
 \begin{center}
 \includegraphics[scale=0.61]{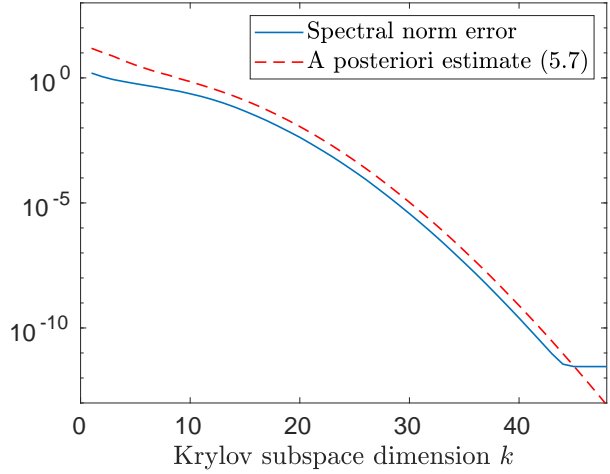}
 \end{center}
 \caption{Convergence of the approximation vs. the a posteriori estimate \eqref{eq:apost_estimate}. }
 \label{fig:bound2}
 \end{figure}

To derive a heuristic a posteriori estimate, we neglect the second term in equation
\eqref{eq:apost_representation} and approximate the first integral by leaving out the exponentials.
This is especially meaningful in the case $A$ has its numerical range on left half-plane, since then the exponentials have their norms less than or equal to 1.
This leads to the approximation
\begin{equation} \label{eq:apost_approx_1}
\begin{aligned}
\mathcal{E}_k(t)  \approx \int_0^t R_k(s) \, \dd s = &U_{k+1} H_{k+1,k} E_k^T \Big(\int_0^t Y_k(s) \, \dd s \Big) V_k^T \\
& +  V_k \Big(\int_0^t Y_k(s) \, \dd s \Big)^T E_k H_{k+1,k}^T U_{k+1}^T.
\end{aligned}
\end{equation}
From a careful inspection we see that $U_{k+1}^T V_k = 0$ implies
\begin{equation} \label{eq:apost_approx_2}
\begin{aligned}
\normm{\int_0^t R_k(s) \, \dd s} & = \normm{ U_{k+1} H_{k+1,k} E_k^T \Big(\int_0^t Y_k(s) \, \dd s \Big) V_k^T} \\
& = \normm{ H_{k+1,k} E_k^T \Big(\int_0^t Y_k(s) \, \dd s \Big) }.
\end{aligned}
\end{equation}
The integral $\int_0^t Y_k(s) \, \dd s$ can be estimated by simply summing
$$
\int_0^t Y_k(s) \, \dd s \approx \sum\limits_{\ell=1}^m \Delta t Y_k(\ell \Delta t),
$$
where $\Delta t = t/m$ and where the intermediate values $Y_k(\ell \Delta t)$ can be 
obtained from the summing and squaring phase
of Algorithm~\ref{Alg:main} (Subsection~\ref{subsec:scaling_and_squaring}).
From \eqref{eq:apost_approx_1} and \eqref{eq:apost_approx_2} we arrive to a computationally efficient a posteriori estimate
\begin{equation} \label{eq:apost_estimate}
est_k(t) := \normm{  H_{k+1,k} E_k^T \sum_{\ell=1}^m \Delta t Y_k(\ell \Delta t) \, }.
\end{equation}
To illustrate the efficiency of this estimate consider the following toy example.
Let $A \in \mathbb{R}^{400 \times 400}$ be the tridiagonal matrix $10^2 \cdot \mathrm{diag}(1,-2,1)$, $t=0.1$, and let
$Z\in \mathbb{R}^{400}$ and $C \in \mathbb{R}^{400}$ be random vectors. Figure \ref{fig:bound2} shows the 
error $\norm{X(t) - X_k(t)}$ vs. the estimate \eqref{eq:apost_estimate}.

We note that by using the error representation \eqref{eq:apost_representation} and the residual $R_k(t)$ given in
\eqref{eq:R_k_representation} it is possible to derive corrected schemes, similarly as is done for the matrix exponential
in~\cite{Botchev} and~\cite{Saad}.

\section{Rank cut}

We describe next a rank cut strategy used in our numerical experiments.
When using Algorithm~\ref{Alg:main}, if
$\mathrm{rank} \begin{bmatrix} X_0 & Q \end{bmatrix} = m$,
then after $k$ iterations the numerical solution 
$X_k(t) = V_k Y_k(t) V_k^T$ has a rank at most $km$ and memory for $\mathcal{O}(kmn)$ entries is needed.
However, the numerical rank of $X_k(t)$ may be considerably smaller
already for small treshold values.
Therefore, when performing multiple time stepping it is reasonable to cut the rank after each step.
This can be done, for example, as follows.
Let $X \in \mathbb{R}^{n \times n}$ and $\sigma_1 \leq \sigma_2 \leq \ldots \leq \sigma_n$ 
denote the singular values of $X$ and $u_i,v_i$ the corresponding left and right singular vectors. 
Consider the projection
$$
P_\epsilon(X) : = \sum_{\sigma_i > \epsilon} \sigma_i u_i v_i^T.
$$
This projection is efficiently applied on the numerical solution $X_k(t)$ given by Algorithm~\ref{Alg:main} since obviously
\begin{equation*} 
P_\varepsilon(X_k(t)) =  V_k P_\varepsilon(Y_k(t)) V_k^T.
\end{equation*}

We have the following bound for the effect of the rank cut of the initial value.
\begin{theorem} \label{thm:cut_svd}
Suppose $\widetilde{X}_0 = P_\varepsilon(X_0)$, and let $X(t)$ and $\widetilde{X}(t)$ be solutions of the system \eqref{eq:DRE}
for initial values $X_0$ and $\widetilde{X}_0$, respectively. Then, 
$$
\norm{X(t) - \widetilde{X}(t) } \leq \varepsilon \, \ee^{2 t \mu(A)}. 
$$
\begin{proof} 
By Thm.~\ref{thm:exact_solution},
\begin{equation} \label{eq:eq_x}
	\begin{aligned}
		X(t) - \widetilde{X}(t) = & \ee^{tA}(X_0 -  \widetilde{X}_0) \ee^{t A^T}  \\
		 & + \int_0^t 
		\ee^{(t-s) A} \big(\widetilde{X}(s) S \widetilde{X}(s)     -   X(s) S X(s)   \big) \ee^{(t-s) A^T} \, \dd s.
	\end{aligned}
\end{equation}
Since $\widetilde{X}_0 = P_\varepsilon(X_0)$, we see that $\widetilde{X}_0 \leq X_0$.
Then, from Thm.~\ref{thm:positivity_and_monotonicity} it follows that $\widetilde{X}(t) \leq X(t)$.
Furthermore, by Lemma~\ref{lem:ABA}, $\widetilde{X}(s) S \widetilde{X}(s) \leq X(s) S X(s)$ for all $s \geq 0$.
Therefore, the second integral on the right hand side of \eqref{eq:eq_x} is a negative semidefinite matrix and thus
$\norm{X(t) - \widetilde{X}(t) } \leq  \norm{ \ee^{tA}(X_0 -  \widetilde{X}_0) \ee^{t A^T} }$, and the claim follows.
\end{proof}
\end{theorem}

A straightforward corollary of Theorem~\ref{thm:cut_svd} is an estimate for the total error arising solely from the rank cutting. 

\begin{corollary}
Consider a time stepping scheme where a rank cut of size $\varepsilon_\ell$ is made after every step $\ell$,
and that the substeps are otherwise exact solutions of \eqref{eq:DRE}.
I.e., the result of step $\ell$ is $X_\ell = P_{\varepsilon_\ell}(X(h))$, where $X(h)$ is the solution of \eqref{eq:DRE} with
the initial value $X_{\ell-1}$. Carrying out this procedure for $n$ steps,
it follows from \ref{thm:cut_svd} and Lady Windermere's fan (see~\cite[Ch.\;I.7]{Hairer}) that the global error satisfies
\begin{equation} \label{eq:svd_estimate}
\norm{ X(n  h) - X_n } \leq \sum\limits_{\ell=1}^n \varepsilon_\ell \, \ee^{2  (n - \ell) h \mu(A) }.
\end{equation}
\end{corollary}

\section{Numerical experiments: optimal cooling problem}

As a numerical example we consider an optimal cooling problem described in~\cite{Saak} (see also Example 2 in~\cite{Stillfjord}).
The underlying linear system is of the form
\begin{equation} \label{eq:linear_steel}
\begin{aligned}
M \dot{x} &= A x + B u, \\
y &= Cx,
\end{aligned}
\end{equation}
where the coefficient matrices arise from a finite element discretization of the cross section of a rail.
A discretization of dimension $n$ gives coefficient matrices $A,M  \in \mathbb{R}^{n \times n}$, $B \in \mathbb{R}^{n \times 7}$ 
and $C \in \mathbb{R}^{n \times 6}$, where $A$ is symmetric. 
This leads to a symmetric DRE of the form \eqref{eq:DRE} with the coefficient matrices
$\widetilde{A} = M^{-1} A$, $Q = C^T C$ and $S = M^{-1}B (M^{-1}B)^T$. We take zero initial value for the DRE. 
The mass matrix $M$ is sparse so the
products using the matrix $M^{-1} A$ are cheap. We note that by a symmetric decomposition of the mass matrix, $M=L^T L$,
the system could also be written
as a system using a symmetric coefficient matrix $L^{-T} A L^{-1} $ for the scaled variable $L X L^{-1}$.

\subsection{Case $n=1357$}

\begin{figure}[h!]
 \begin{center}
 \includegraphics[scale=0.61]{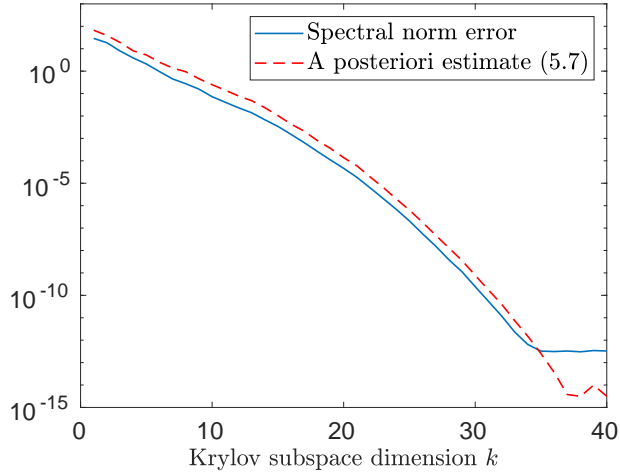}
 \end{center}
 \caption{Convergence of the approximation given by Algorithm~\ref{Alg:main} and its a posteriori estimate \eqref{eq:apost_estimate}.}
 \label{fig:rail1}
 \end{figure}

Figure~\ref{fig:rail1} shows the convergence of Algorithm~\ref{Alg:main} and an a posteriori error estimate 
given by \eqref{eq:apost_estimate}, when $T=10$. We compute the spectral norm error $\norm{X(T) - X_k(T)}$
for different Krylov subspace dimensions $k$.
For the scaling and squaring part (Subsection~\ref{subsec:scaling_and_squaring}), we set the parameter $m=10$.
Table~\ref{table:ex1} shows the CPU time needed for the block Krylov process and for the scaling squaring part of Algorithm~\ref{Alg:main},
for four different Krylov subspace sizes.

\begin{figure}[h!]
 \begin{center}
 \includegraphics[scale=0.61]{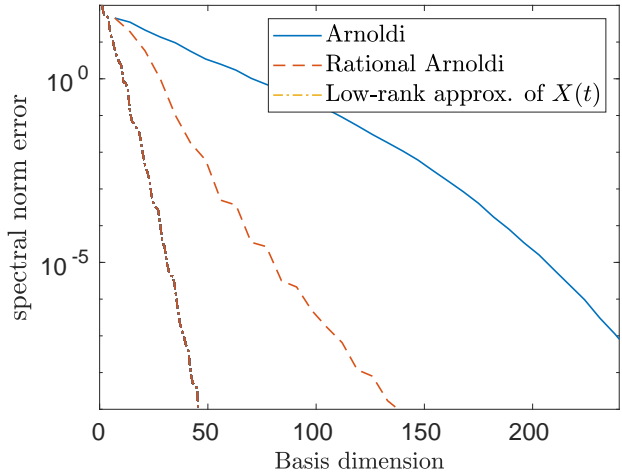}
 \end{center}
 \caption{Convergence of a single step approximation with the Arnoldi and the rational Arnoldi iteration. The figure shows also the convergence 
 of the best low rank approximation of $X(t)$. }
 \label{fig:rail_rational}
 \end{figure}

Figure~\ref{fig:rail_rational} shows the convergence of a single step for $T=20$, when we apply the block orthogonalisation
procedure of Subsection~\ref{subsec:Krylov_approx} on the Krylov subspace \eqref{def:krylov} and on a rational Krylov subspace \eqref{def:rational}
spanned by $A$. For the rational Krylov subspace we set all nodes $s_i$ equal to 1. 
Here subspace dimension denotes the number of columns of the basis matrix $V_k$.
For comparison, we also consider the best
low rank approximation of the solution $X(T)$ obtained from its singular value decomposition (SVD) for different
ranks (denoted basis dimension in the figure). We see that the rational
approximation needs a considerably smaller subspace for a given error than the polynomial approximation.

\begin{table}[h!] 
\begin{center}
\caption{Timings for the Krylov subspace iteration and for the solving of the projected system using the modified
 Davison--Maki method, when integrating up to $t=10$ using a single Krylov subspace iteration. Times are in seconds.}
\begin{tabular}{ccc}
 $k$ & Krylov iteration & solving small dimensional system  \\
 \hline
  $10$  &   0.046  &  0.037    \\
  $20$  &   0.16   &  0.091    \\
  $30$  &   0.31 &  0.20    \\
  $40$  &   0.49 &  0.44    \\
 \end{tabular}
 \label{table:ex1}
\end{center}
\end{table}

\begin{figure}[h!]
 \begin{center}
 \includegraphics[scale=0.61]{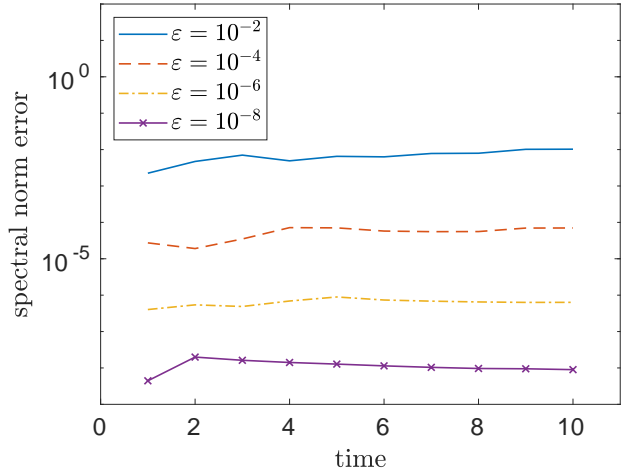}
 \end{center}
 \caption{The error of the numerical solution for different tolerances $\varepsilon$.}
 \label{fig:rail2}
 \end{figure}

\begin{figure}[h!]
 \begin{center}
 \includegraphics[scale=0.61]{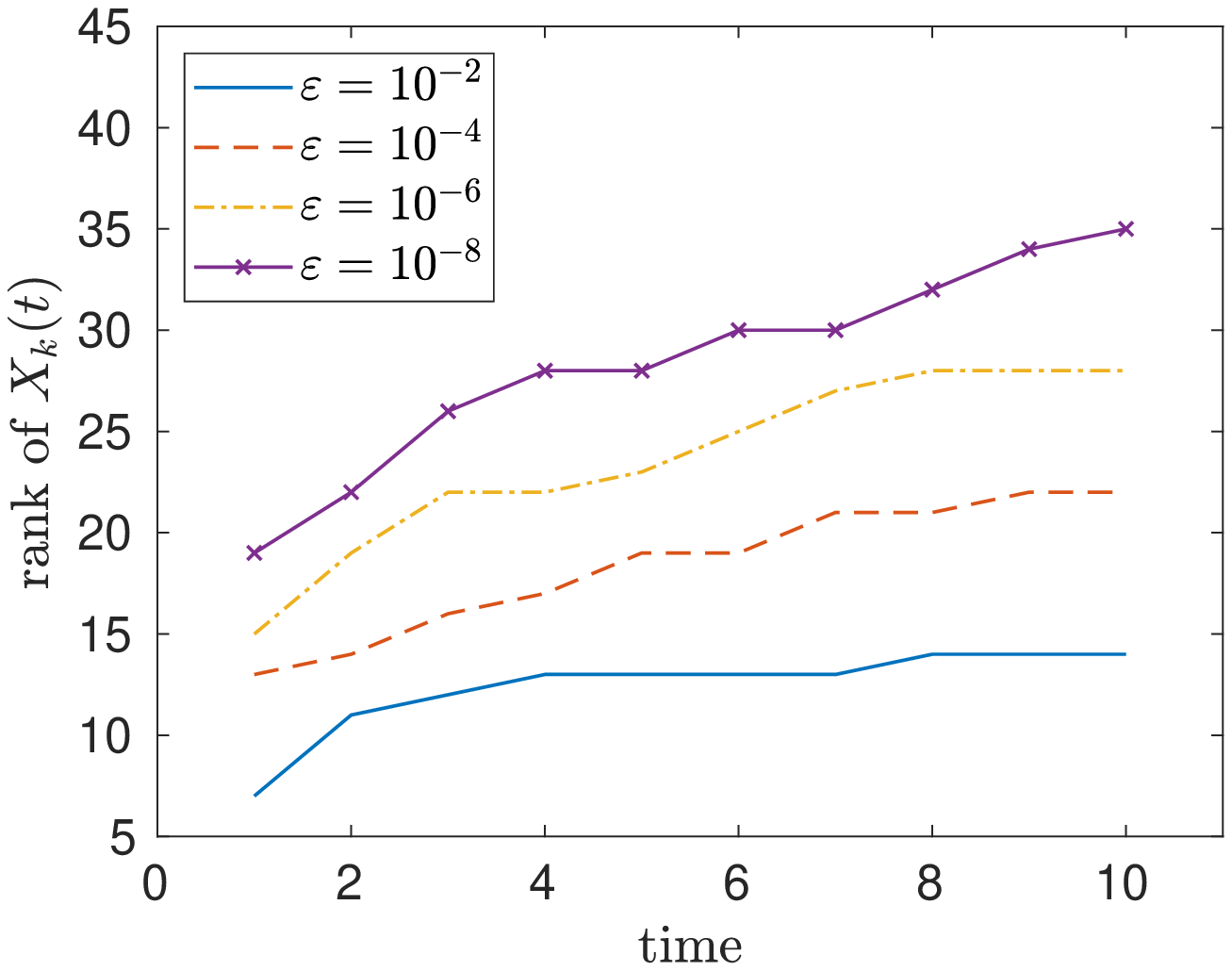}
 \end{center}
 \caption{The rank of the numerical solution for different tolerances $\varepsilon$.}
 \label{fig:rail3}
 \end{figure}

Next, we apply Algorithm~\ref{Alg:main} for $N=10$ subsequent steps. We set for the Krylov error a tolerance $\varepsilon$,
and use the a posteriori estimate \eqref{eq:apost_estimate} as a criterion for stopping the iteration.
Also, after each step we cut the rank using the projector $P_\varepsilon$.
Figure~\ref{fig:rail2} depicts the final errors at $T=10$ for 4 different values of $\varepsilon$. As we see the
final errors are not far from the tolerances $\varepsilon$ used for substeps. Figure~\ref{fig:rail3} depicts the
growth of the rank in the numerical solution for different tolerances $\varepsilon$. 
We see that the substepping approach requires less memory
for a given error tolerance than a single run using Algorithm~\ref{Alg:main}. This is depicted in Table~\ref{table:ex2}.

\begin{table}[h!] 
\begin{center}
\caption{Maximum number of columns of the basis matrix $V_k$ along the iteration for the substepping approach 
 and for one step approximation using Algorithm~\ref{Alg:main}, when an error tolerance $\mathrm{tol}$ is required.}
\begin{tabular}{ccc}
 $\mathrm{tol}$ &  time stepping &  single step iteration  \\
 \hline
  $10^{-2}$  &   60  &  112   \\
  $10^{-4}$  &   76  &  147   \\
  $10^{-6}$  &   112 &  175   \\
  $10^{-8}$  &   160 &  203   \\
 \end{tabular} 
 \label{table:ex2}
\end{center}
\end{table}

\begin{figure}[h!]
 \begin{center}
 \includegraphics[scale=0.61]{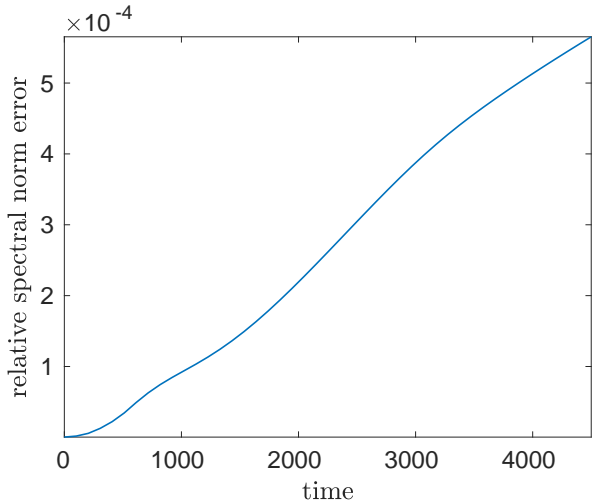}
 \end{center}
 \caption{The relative spectral norm error of the numerical solution in a time integration up to $T=4500$.}
 \label{fig:rail4}
 \end{figure}

As a last experiment for the case $n=1357$, we carry out a time integration up to $T=4500$ using $900$ substeps.
Figure~\ref{fig:rail4} shows the relative spectral norm error along the time integration,
i.e., the error $\norm{\widetilde{X}(t) - X(t)}/\norm{X(t)}$, where $\widetilde{X}(t)$ denotes the numerical solution.
We use a Krylov subspace dimension $k=32$ for the first substep and $k=20$ for the rest.
After each time step, we cut the rank of the numerical solution to $40$ using SVD. By these choices of Krylov subspace sizes we ensure that
the error arising from the rank cut dominates at each time step.

\begin{figure}[h!]
\begin{center}
\includegraphics[scale=0.61]{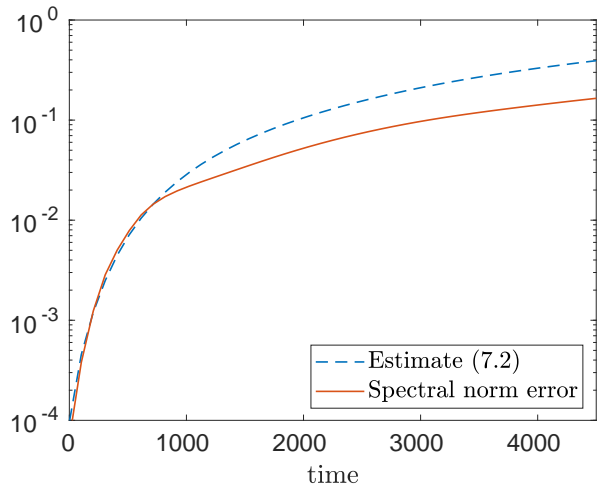}
\end{center}
\caption{The spectral norm error of the approximation and the estimate \eqref{eq:est_num} for time integration up to $T=4500$.}
\label{fig:rail4_svd}
\end{figure}

Next, in order to use the estimate \eqref{eq:svd_estimate}, we approximate $\mu(A) \approx 0$ and assuming that the error arising from the rank cut dominates the total error.
We then approximate the total using \eqref{eq:svd_estimate} as
\begin{equation} \label{eq:est_num}
\norm{ X(n  h) - X_n } \approx \sum\limits_{\ell=1}^n \varepsilon_\ell.
\end{equation}
Figure~\ref{fig:rail4_svd} shows the error arising from the best 2-norm approximation after each step, i.e., the singular value $\sigma_{41}$,
and the estimate \eqref{eq:est_num}. We see that the error in the end is not far from $900 \cdot \sigma_{41}^{(900)}$, the number of time steps times the largest
rank cut.

\subsection{Case $n=5177$}

\begin{figure}[h!]
\begin{center}
\includegraphics[scale=0.61]{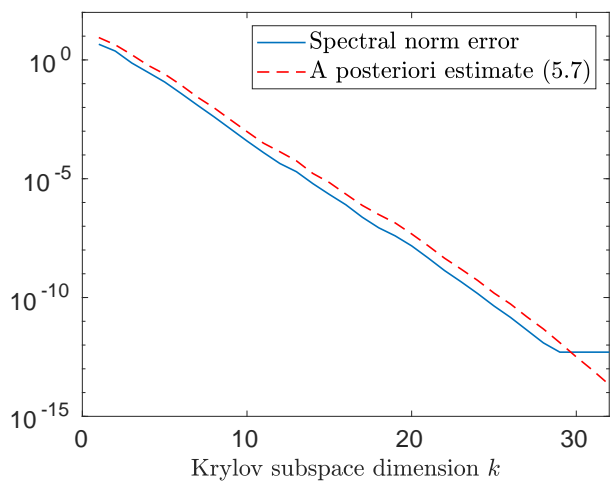}
\end{center}
 \caption{Convergence of the approximation given by Algorithm~\ref{Alg:main} and its a posteriori estimate \eqref{eq:apost_estimate}, when $n=5177$. }
\label{fig:rail5}
\end{figure}

Next we consider a finite element discretization with $n=5177$. 
Figure~\ref{fig:rail5} shows the convergence of Algorithm~\ref{Alg:main} and an a posteriori error estimate 
given by \eqref{eq:apost_estimate}, when $T=5$.
For the scaling and squaring part we set the parameter $m=10$.

\begin{figure}[h!]
\begin{center}
\includegraphics[scale=0.61]{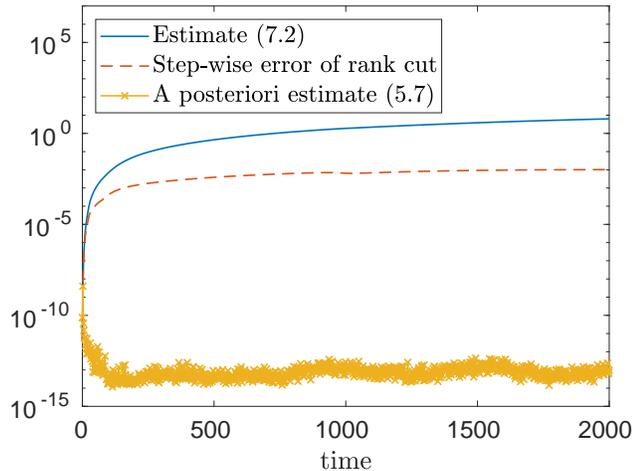}
\end{center}
\caption{The estimate \eqref{eq:est_num} of the spectral norm error, size of the rank cut at each time step 
and the a posteriori estimate \eqref{eq:apost_estimate} at each time step, when $n=5177$ and time integration up to $T=2000$.}
\label{fig:convergence_5177}
\end{figure}

We next carry out a time integration up to $T=2000$ using $1000$ substeps. We estimate the total error without access to a reference solution
using the estimate \eqref{eq:est_num}.
As above, we use a Krylov subspace dimension $k=32$ for the first step, and $k=20$ for rest of the steps,
and cut the rank to 40 after each step using SVD. We see from Figure~\ref{fig:convergence_5177} that the a posteriori error
estimate for Algorithm~\ref{Alg:main} is negligible compared to the error arising from the best 2-norm approximation at each step.
Figure~\ref{fig:convergence_5177} shows also the estimate \eqref{eq:est_num}. 
We see that the estimate is of the same order ($\approx$ 10 times bigger) as in the $n=1357$-case.

\section{Conclusions and Outlook}

We have proposed a Krylov subspace approximation method for large scale differential Riccati equations. 
We have proven that the method is structure preserving in the sense that it preserves 
two important properties of the exact flow, namely the property of the positivity and also under certain 
conditions also the property of the monotonicity. 
We have also provided an a priori error analysis of the Krylov subspace approximation which 
shows a superlinear convergence. This behavior was also verified in numerical experiments.  
In addition, an a posteriori error analysis was carried out and the proposed estimate was shown to be accurate in numerical examples.
In order to limit the memory consumption, we considered limiting the rank of the numerical solution in multiple time stepping.
To avoid excessively large approximation basis $V_k$, more studies of the rational
Krylov subspaces are needed. Their benefits were illustrated in numerical experiments.

We would like to point out that the presented block Krylov subspace method can be extended to the unsymmetric differential Riccati equation. 
A possible extension could also be the nonautonomous case, i.e, the case in which the coefficient matrices are $Q$, $S$ and $A$ are time dependent.
In this case an essential tool would be the so called Magnus expansion (see e.g.~\cite{Bader}) which gives the fundamental solution of the linear system
corresponding to the time dependent coefficient matrix $A$.

\section*{Acknowledgments}
The authors thank Valeria Simoncini for pointing out relevant literature related to the algebraic Riccati equation
and Tony Stillfjord for several helpful comments on a draft of the paper.

\appendix

\section{Auxiliary Lemmas and the proof of Thm.~\ref{thm:main_conv}} 

\vspace{5mm}

We first state two lemmas needed in Thm.~\ref{thm:cut_svd} and Thm.~\ref{thm:main_conv}, respectively. \\


\begin{lemma} \label{lem:ABA}
Let $A,\widetilde{A},B$ be symmetric positive semidefinite matrices
such that $A \leq \widetilde{A}$. Then, also
$$
ABA \leq \widetilde{A} B  \widetilde{A}.
$$
\begin{proof}
Assume first that $B$ is positive definite.
Then, (see~\cite[p.\;431]{Horn})
$$
B^{\frac{1}{2}}  A B^{\frac{1}{2}} \leq B^{\frac{1}{2}} \widetilde{A} B^{\frac{1}{2}}
$$
and therefore also (see~\cite[p.\;438]{Horn})
$$
\Big(B^{\frac{1}{2}} A B^{\frac{1}{2}} \Big)^2 \leq \Big( B^{\frac{1}{2}} \widetilde{A} B^{\frac{1}{2}} \Big)^2.
$$
Then, we see that
\begin{equation*}
	\begin{aligned}
		0 \, \leq &  \, B^{-\frac{1}{2}} \left[
		\Big( B^{\frac{1}{2}}  \widetilde{A} B^{\frac{1}{2}} \Big)^2 -
		\Big( B^{\frac{1}{2}} A B^{\frac{1}{2}} \Big)^2 \right] B^{-\frac{1}{2}} \\
		 = & \, \widetilde{A} B  \widetilde{A} - ABA.
	\end{aligned}
\end{equation*} 
Then, consider the matrix $B_\varepsilon = B + \varepsilon I$, $\varepsilon > 0$, where $B$ is positive semidefinite. Clearly, $B_\varepsilon$ is positive definite for all
$\varepsilon > 0$. Therefore 
\begin{equation} \label{eq:Aeps}
	A B_\varepsilon A \leq \widetilde{A}  B_\varepsilon \widetilde{A} \quad \textrm{for all} \quad \varepsilon > 0.
\end{equation}
Taking the limit $\varepsilon \rightarrow 0$ and using the fact that
$$
\max\limits_{\norm{x} = 1} \, x^T \left( \widetilde{A} B_\varepsilon  \widetilde{A} - A B_\varepsilon A \right) x
$$
is a continuous function of $\varepsilon$, the claim follows.
\end{proof}
\end{lemma}

\begin{lemma} \label{lem:auxiliary1}
Let $A\in \mathbb{R}^{n \times n}$, $B \in \mathbb{R}^{n \times \ell}$,
let $V_k$ be a matrix with orthonormal columns such that $\mathcal{K}_k(A,B) \subset \mathrm{R}(V_k)$ and let $H_k = V_k^T A V_k$.
Then, for all $t,s>0$ it holds that
$$
\norm{ \left(\ee^{tA} - V_k \ee^{t H_k} V_k^T \right) \ee^{s A} B } \leq 4 \, \max\{1,\ee^{(t+s) \mu(A)} \} \frac{\norm{(t+s)A}^k}{k !} \norm{B}.
$$
\begin{proof} Using the polynomial approximation property of the Krylov approximation (see~\cite[Lemma\;3.1]{Saad}), we see that
\begin{equation*}
\begin{aligned}
V_k V_k^T \ee^{s A} B & = V_k V_k^T \sum\limits_{\ell=0}^{k-1} \frac{(sA)^\ell}{\ell !} B
+ V_k V_k^T \sum\limits_{\ell=k}^\infty \frac{(sA)^\ell}{\ell !} B\\
& =  V_k \sum\limits_{\ell=0}^{k-1} \frac{(s H_k)^\ell}{\ell !} V_k^T B + V_k V_k^T \sum\limits_{\ell=k}^\infty \frac{(sA)^\ell}{\ell !}B \\
& = V_k \ee^{s H_k} V_k^T B - V_k \sum\limits_{\ell=k}^\infty \frac{(s H_k)^\ell}{\ell !} V_k^T B
+ V_k V_k^T \sum\limits_{\ell=k}^\infty \frac{(sA)^\ell}{\ell !} B.
\end{aligned}
\end{equation*}
Therefore,
\begin{equation} \label{eq:4_terms}
\begin{aligned}
\left(\ee^{tA} - V_k \ee^{t H_k} V_k^T \right) \ee^{s A} B &= \ee^{(t+s)A} B -  V_k \ee^{t H_k} V_k^T ( V_k V_k^T \ee^{s A} B ) \\
&= \ee^{(t+s)A} B -  V_k \ee^{(t+s) H_k} V_k^T B  \\
 & + V_k \ee^{t H_k} V_k^T  \sum\limits_{\ell=k}^\infty \frac{(sA)^\ell}{\ell !} B
 - V_k \ee^{t H_k} \sum\limits_{\ell=k}^\infty \frac{(s H_k)^\ell}{\ell !} V_k^T B \\
&=\sum\limits_{\ell=k}^\infty \frac{ (\big(t+s \big)A)^\ell }{\ell !} B -  
V_k\sum\limits_{\ell=k}^\infty \frac{ \big((t+s) H_k \big)^\ell }{\ell !} V_k^T B   \\
& + V_k \ee^{t H_k} V_k^T  \sum\limits_{\ell=k}^\infty \frac{(sA)^\ell}{\ell !} B
- V_k \ee^{t H_k} \sum\limits_{\ell=k}^\infty \frac{(s H_k)^\ell}{\ell !} V_k^T B.
\end{aligned}
\end{equation}
Using the bounds $\norm{\ee^{t A}} \leq \ee^{t \mu(A)}$, $\mu(H_k) \leq \mu(A)$ and (see~\cite[Lemma\;A.2]{GallopoulosSaad})
$$
\normm{\sum\limits_{\ell=k}^\infty \frac{(t A)^\ell}{\ell !}} 
\leq \max\{1,\ee^{t \mu(A)} \} \frac{ \norm{tA}^k}{k!} \quad \textrm{for all} \quad t\geq 0
$$
on the four terms on the RHS of \eqref{eq:4_terms}, the claim follows. \end{proof}
\end{lemma}

\begin{lemma} \label{lem:auxiliary2}
Let $X(s)$ be the solution of the Riccati differential equation \eqref{eq:DRE} at time $s$,  $0\leq s \leq t$, and let $V_k$ be a matrix with orthonormal columns
such that $\mathcal{K}_k(A, \begin{bmatrix} C & Z \end{bmatrix})  \subset \mathrm{R}(V_k)$. Denote $H_k = V_k^T A V_k$. Then, the following bound holds:
\begin{equation*}
\begin{aligned}
 & \norm{ \big( \ee^{(t-s)A} - V_k \ee^{(t-s) H_k} V_k^T \big) X(s) } \\
\leq \,\,\, & 4 \, c(s)
 \max \{1, \ee^{(t+s) \mu(A)} \} \norm{A}^k
\left( \frac{t^k}{k!}\norm{X_0} + \frac{t^{k+1}}{(k+1)!} \norm{Q} \right),
\end{aligned}
\end{equation*}
where
$$
c(s) := 1 +  s \norm{S} \max_{w \in [0,s]} \norm{X(w)} \, \varphi_1 \left( s \norm{S} \max_{w \in [0,s]} \norm{X(w)} \, \max \{1, \ee^{ s  \mu(A)} \} \right).
$$
\begin{proof} Using the integral representation \eqref{eq:exact_solution} for $X(s)$ we may write
\begin{equation} \label{eq:split_C_1_C_2} 
\big( \ee^{(t-s)A} - V_k \ee^{(t-s) H_k} V_k^T \big) X(s)  =  c_{1,k}(t,s) + c_{2,k}(t,s),
\end{equation}
where
\begin{equation} \label{eq:C_1}
\begin{aligned}
  c_{1,k}(t,s) &= \left[ \big( \ee^{(t-s) A} - V_k \ee^{(t-s)H_k} V_k^T \big) \ee^{sA} Z \right] Z^T \ee^{s A^T}  \\
& \quad + \int\limits_0^s \left[ \big( \ee^{(t-s)A} - V_k \ee^{(t-s)H_k} V_k^T \big) 
\ee^{(s-u) A} C \right] C^T \ee^{(s-u)A^T} \, \dd u
\end{aligned}
\end{equation}
and
\begin{equation} \label{eq:C_2}
\begin{aligned}
  c_{2,k}(t,s) &= \big( \ee^{(t-s)A} - V_k \ee^{(t-s) H_k} V_k^T \big) 
  \int\limits_0^s \ee^{(s-u)A } X(u) S X(u) \ee^{ (s-u) A^T} \, \dd u.
\end{aligned}
\end{equation} 
By using Lemma~\ref{lem:auxiliary1}, we obtain for the expressions inside the square brackets on
right hand side of \eqref{eq:C_1} the bounds
\begin{equation}
\begin{aligned}
&\normm{ \left[ \big( \ee^{(t-s) A} - V_k \ee^{(t-s)H_k} V_k^T \big) \ee^{sA} Z \right] Z^T \ee^{s A^T} } \\
 \leq & \,\,\,  4  \max \{1,\ee^{(t+s) \mu(A)}\} \norm{X_0} \norm{A}^k \frac{t^k}{k!}
\end{aligned}
\end{equation}
and
\begin{equation*}
\begin{aligned}
& \normm{\int\limits_0^s \left[ \big( \ee^{(t-s)A} - V_k \ee^{(t-s)H_k} V_k^T \big) 
\ee^{(s-u) A} C \right] C^T \ee^{(s-u)A^T} \, \dd u } \\ 
 \leq &\,\,\,  4  \norm{Q} \int\limits_0^s \frac{\big((t-u)\norm{A} \big)^k}{k!} 
\max \{1, \ee^{(t-u) \mu(A)} \} \max \{1, \ee^{(s-u) \mu(A)} \} \, \dd u \\
\leq &\,\,\,  4 \norm{Q}  \max \{1, \ee^{(t+s) \mu(A)} \} \norm{A}^k \frac{ t^{k+1}}{(k+1)!} .
\end{aligned}
\end{equation*}
Thus,
\begin{equation} \label{eq:C_1_bound}
\normm{c_{1,k}(t,s)} \leq 4 \max \{1, \ee^{(t+s) \mu(A)} \} \norm{A}^k \left( \frac{t^k}{k!} \norm{X_0} + \frac{t^{k+1}}{(k+1)!} \norm{Q} \right).
\end{equation}
From \eqref{eq:C_2} we see that
\begin{equation} \label{eq:C_2_first_step}
\begin{aligned}
 \normm{c_{2,k}(t,s)} &\leq \int\limits_0^s \normm{ \big( \ee^{(t-s)A} - V_k \ee^{(t-s) H_k} V_k^T \big) \ee^{(s-u)A } X(u) } \\
  & \cdot \normm{S} \max_{w \in [0,s]} \norm{X(w)} \ee^{(s-u)\mu(A)} \, \dd u.
 \end{aligned}
\end{equation}
Next we bound the first factor in the integrand of \eqref{eq:C_2_first_step}.
We substitute the integral representation \eqref{eq:exact_solution} for $X(u)$ to find that
\begin{equation} \label{eq:C_2_second_step}
\begin{aligned}
&\normm{( \ee^{(t-s)A} - V_k \ee^{(t-s) H_k} V_k^T ) \ee^{(s-u)A} X(u)} \\
\leq & \normm{ \left[ ( \ee^{(t-s)A} - V_k \ee^{(t-s) H_k} V_k^T ) \ee^{sA} Z \right] Z^T \ee^{u A^T}} \\
+ & \int\limits_0^u \normm{\left[ ( \ee^{(t-s)A} - V_k \ee^{(t-s) H_k} V_k^T ) \ee^{(s-w)A} C \right] C^T \ee^{(u-w)A^T} } \, \dd w\\
 + & \int\limits_0^u \normm{ ( \ee^{(t-s)A} - V_k \ee^{(t-s) H_k} V_k^T ) \ee^{(s-w)A} X(w) } \\
 & \cdot \norm{S} \max_{w \in [0,u]} \norm{X(w)} \max \{1, \ee^{(u-w) \mu(A)} \} \, \dd w.
\end{aligned}
\end{equation}
As above when bounding $\norm{c_{1,k}(t,s)}$, we use Lemma~\ref{lem:auxiliary1} on the expressions inside the square brackets on
right hand side of \eqref{eq:C_2_second_step}, to get the inequality
\begin{equation} \label{eq:C_2_third_step}
\begin{aligned}
&\normm{( \ee^{(t-s)A} - V_k \ee^{(t-s) H_k} V_k^T ) \ee^{(s-u)A} X(u)} \\
\leq &  \,\, 4 \norm{A}^k \max \{1, \ee^{(t+u) \mu(A)} \} \left( \frac{t^k}{k!}\norm{X_0} + \frac{t^{k+1}}{(k+1)!} \norm{Q} \right) \\
 & +  \norm{S} \max_{w \in [0,u]} \norm{X(w)} \, \max \{1, \ee^{u \mu(A)} \} \\
 & \cdot \int\limits_0^u \normm{ ( \ee^{(t-s)A} - V_k \ee^{(t-s) H_k} V_k^T ) \ee^{(s-w)A} X(w) }  \, \dd w.
\end{aligned}
\end{equation}
Applying Gr\"onwall's lemma on \eqref{eq:C_2_third_step}, we find that
\begin{equation} \label{eq:C_2_fourth_step}
\begin{aligned}
&\normm{( \ee^{(t-s)A} - V_k \ee^{(t-s) H_k} V_k^T ) \ee^{(s-u)A} X(u)} \\
\leq &  \,\, 4 \norm{A}^k \max \{1, \ee^{(t+u) \mu(A)} \} \left( \frac{t^k}{k!}\norm{X_0} + \frac{t^{k+1}}{(k+1)!} \norm{Q} \right) \\
 & \cdot \ee^{ u \norm{S} \max_{w \in [0,u]} \norm{X(w)} \, \max \{1, \ee^{u \mu(A)} \}}.
\end{aligned}
\end{equation}
Substituting \eqref{eq:C_2_fourth_step} into \eqref{eq:C_2_first_step}, we get
\begin{equation} \label{eq:C_2_bound}
\begin{aligned}
& \normm{c_{2,k}(t,s)} \leq  4 \norm{S} \max_{w \in [0,s]} \norm{X(w)}  \norm{A}^k \max \{1, \ee^{(t+s) \mu(A)} \} \\
& \cdot \left( \frac{t^k}{k!}\norm{X_0} + \frac{t^{k+1}}{(k+1)!} \norm{Q} \right)  
s \varphi_1 \left( s \norm{S} \max_{w \in [0,s]} \norm{X(w)} \, \max \{1, \ee^{ s  \mu(A)} \} \right).
\end{aligned}
\end{equation}
The bounds \eqref{eq:C_1_bound} and \eqref{eq:C_2_bound} together show the claim. \end{proof} 
\end{lemma}

Using Lemmas~\ref{lem:auxiliary1} and~\ref{lem:auxiliary2} we are now ready to prove Theorem~\ref{thm:main_conv}.

\subsection*{Proof of Theorem~\ref{thm:main_conv}}

\begin{proof} From the integral representation \eqref{eq:exact_solution} for $X(t)$
and for the solution $Y_k(t)$ of the small dimensional system \eqref{eq:riccati_small_system}, 
we see that
\begin{equation} \label{eq:error_div}
X(t) - X_k(t) = F_{1,k}(t) + F_{2,k}(t),
\end{equation}
where
\begin{equation*}
\begin{aligned}
F_{1,k}(t) := & \,\,\, \ee^{t A} X_0 \ee^{t A^T} - V_k \ee^{ t H_k}  V_k^T X_0  V_k \ee^{ t H_k^T} V_k^T   \\
& + \quad \int\limits_0^t \Big( \ee^{(t-s) A} Q \ee^{(t-s) A^T} - V_k \ee^{ (t-s) H_k} Q_k \ee^{ (t-s) H_k^T} V_k^T \Big)\, \dd s,
\end{aligned}
\end{equation*}
and
\begin{equation} \label{eq:F_2_k}
 \begin{aligned}
F_{2,k}(t) &= \int\limits_0^t \ee^{(t-s) A} X(s) S X(s) \ee^{(t-s) A^T} \, \dd s \\
& -  \int\limits_0^t V_k \ee^{(t-s) H_k} V_k^T X_k(s) S X_k(s) V_k \ee^{(t-s) H_k^T} V_k^T \, \dd s.  
 \end{aligned}
\end{equation}
Theorem~\ref{thm:Lyapunov} shows that $F_{1,k}(t)$ is bounded as
\begin{equation} \label{eq:F_1_bound}
 \begin{aligned}
\norm{ F_{1,k}(t)} \leq  4  \max \{1, \ee^{2 t \mu(A)} \} \norm{A}^k \left( \frac{t^k}{k!}\norm{X_0} + \frac{t^{k+1}}{(k+1)!} \norm{Q} \right).
 \end{aligned}
\end{equation}
We add and substract the term
$$
\int\limits_0^t \ee^{(t-s) A} X(s) S  X_k(s) V_k \ee^{(t-s) H_k^T} V_k^T \, \dd s
$$ 
to \eqref{eq:F_2_k} to obtain
\begin{equation} \label{eq:F_2_k_2}
 \begin{aligned}
F_{2,k}(t) &= \int\limits_0^t \ee^{(t-s) A} X(s) S \, F_{3,k}(t,s)^T \, \dd s   \\
 & + \int\limits_0^t F_{3,k}(t,s) \, S \,X_k(s) V_k \ee^{(t-s) H_k^T} V_k^T \, \dd s,
 \end{aligned}
\end{equation}
where
\begin{equation} \label{eq:F_3_k}
\begin{aligned}
& F_{3,k}(t,s) = \ee^{(t-s) A } X(s)  -   V_k \ee^{(t-s) H_k} V_k^T  X_k(s) \\
&= \big(\ee^{(t-s) A} - V_k \ee^{(t-s) H_k} V_k^T \big) X(s) - V_k \ee^{(t-s) H_k} V_k^T \big( X_k(s) - X(s) \big).
\end{aligned}
\end{equation}
From \eqref{eq:F_2_k_2} and \eqref{eq:F_3_k} we see that
\begin{equation} \label{eq:F_2_bound}
 \begin{aligned}
 \norm{ F_{2,k}(t) } \leq & \,  2  \, \norm{S} \, \alpha(t) \int\limits_0^t \max \{1, \ee^{(t-s) \mu(A)} \}  \\ 
& \cdot \Big( \norm{\big(\ee^{(t-s) A} - V_k \ee^{(t-s) H_k} V_k^T \big) X(s)}  + \norm{X(s) - X_k(s)} \Big) \, \dd s.
 \end{aligned}
\end{equation}
where 
$$
\alpha(t) = \max \left\{ \max_{s \in [0,t]} \norm{X(s)},\max_{s \in [0,t]} \norm{X_k(s)} \right\}.
$$
The claim follows now from \eqref{eq:error_div}, \eqref{eq:F_1_bound}, \eqref{eq:F_2_bound}, Lemma~\ref{lem:auxiliary1}, Gr\"onwall's lemma,
Corollary~\ref{cor:max_X} and Corollary~\ref{cor:max_X_k}, which form a sequence of substitutions.\end{proof}



 \end{document}